\documentclass[12pt]{article}
\usepackage{graphicx}
\usepackage{color}
\catcode`\@=11 \@addtoreset{equation}{section}

\catcode`\@=12

\newtheorem{Theorem}{Theorem}[section]
\newtheorem{Proposition}{Proposition}[section]
\newtheorem{Lemma}{Lemma}[section]
\newtheorem{Corollary}{Corollary}[section]
\newtheorem{Remark}{Remark}[section]

\newcommand{\bTheorem}[1]{
\begin{Theorem} \label{T#1} }
\newcommand{\eT}{\end{Theorem}}

\newcommand{\bProposition}[1]{
\begin{Proposition} \label{P#1}}
\newcommand{\eP}{\end{Proposition}}

\newcommand{\bLemma}[1]{
\begin{Lemma} \label{L#1} }
\newcommand{\eL}{\end{Lemma}}

\newcommand{\bCorollary}[1]{
\begin{Corollary} \label{C#1} }
\newcommand{\eC}{\end{Corollary}}

\newcommand{\bRemark}[1]{
\begin{Remark} \label{R#1} }
\newcommand{\eR}{\end{Remark}}

\newcommand{\bFormula}[1]{
\begin{equation} \label{#1}}
\newcommand{\eF}{\end{equation}}

\newcommand{\Ov}[1]{\overline{#1}}

\newcommand{\DC}{C^\infty_c}

\newcommand{\vr}{\varrho}

\newcommand{\vu}{\vc{u}}
\newcommand{\vc}[1]{{\bf #1}}
\newcommand{\qed}{\bigskip \rightline {Q.E.D.} \bigskip}
\newcommand{\Div}{{\rm div}_x}
\newcommand{\Grad}{\nabla_x}

\newcommand{\tn}[1]{\mbox {\F #1}}
\newcommand{\dx}{{\rm d} {x}}
\newcommand{\dt}{{\rm d} t }

\newcommand{\intO}[1]{\int_{\Omega} #1 \ \dx}

\newcommand{\bProof}{{\bf Proof: }}

\newcommand{\ep}{\varepsilon}

\font\F=msbm10 scaled 1000

\definecolor{Cgrey}{rgb}{0.85,0.85,0.85}
\definecolor{Cblue}{rgb}{0.50,0.85,0.85}
\definecolor{Cred}{rgb}{1,0,0}
\definecolor{fancy}{rgb}{0.10,0.85,0.10}

\newcommand\Cbox[2]{%
    \newbox\contentbox%
    \newbox\bkgdbox%
    \setbox\contentbox\hbox to \hsize{%
        \vtop{
            \kern\columnsep
            \hbox to \hsize{%
                \kern\columnsep%
                \advance\hsize by -2\columnsep%
                \setlength{\textwidth}{\hsize}%
                \vbox{
                    \parskip=\baselineskip
                    \parindent=0bp
                    #2
                }%
                \kern\columnsep%
            }%
            \kern\columnsep%
        }%
    }%
    \setbox\bkgdbox\vbox{
        \color{#1}
        \hrule width  \wd\contentbox %
               height \ht\contentbox %
               depth  \dp\contentbox
        \color{black}
    }%
    \wd\bkgdbox=0bp%
    \vbox{\hbox to \hsize{\box\bkgdbox\box\contentbox}}%
    \vskip\baselineskip%
}

\date{}

\begin{document}


\title{Maximal dissipation and well-posedness for the compressible Euler system}

\author{Eduard Feireisl \thanks{The research of E.F. leading to these results has received funding from the European Research Council under the European Union's Seventh Framework
Programme (FP7/2007-2013)/ ERC Grant Agreement 320078.} }

\maketitle

\bigskip

\centerline{Institute of Mathematics of the Academy of Sciences of the Czech Republic}

\centerline{\v Zitn\' a 25, 115 67 Praha 1, Czech Republic}






\maketitle

\bigskip





\begin{abstract}

We discuss the problem of well-posedness of the compressible (barotropic) Euler system in the framework of weak solutions. The principle of
\emph{maximal dissipation} introduced by C.M. Dafermos is adapted and combined with the concept of \emph{admissible} weak solutions. We use the method
of convex integration in the spirit of the recent work of C.DeLellis and L.Sz\' ekelyhidi to show various counterexamples to well-posedness. On the other hand, we conjecture that the principle of \emph{maximal dissipation} should be retained as a possible criterion of uniqueness as it is violated by the oscillatory solutions obtained in the process of convex integration.

\end{abstract}

\medskip


\bigskip


\section{Introduction}
\label{i}

The problem of \emph{well-posedness} for general  systems of  (nonlinear) hyperbolic equations remains largely open, despite an enormous amount of literature
and a few particular situations, where rigorous results are available, see e.g. Benzoni-Gavage, Serre
\cite{BenSer}, Bressan \cite{BRESSAN}, Dafermos \cite{D4}, LeFloch \cite{LeFl}, Serre \cite{SERRE1} for a
review of the state-of-art. The inevitable presence of singularities that may develop in a finite time, no matter how smooth and/or ``small'' the data are, gave rise to several concepts of weak (distributional) solutions, supplemented with various admissibility criteria to pick up the physically relevant solution.
As is well known, nonlinear hyperbolic systems are not well-posed in the class of weak solutions and
may admit typically infinitely many distributional solutions emanating from the same initial data, among them some
apparently non-physical violating certain underlying principles as, for instance, the Second law of thermodynamics.

The issue of admissibility of weak solutions has been examined recently in the light of the new results of DeLellis and Sz\' ekelyhidy \cite{DelSze3}, \cite{DelSze}
in the context of  gas dynamics. Adapting the method of \emph{convex integration} (see M\" uller and \v Sver\' ak \cite{MullSve}) DeLellis an Sz\' ekelyhidy constructed infinitely many solutions of the incompressible Euler system emanating from the same initial data and satisfying the standard admissibility criterion based on mechanical energy dissipation. These results provided simple examples of non-uniqueness even in the context of compressible fluids and were later
extended in this direction by Chiodaroli \cite{Chiod}.

The weak solutions obtained via convex integration are ``oscillatory'', or more precisely, their construction is based on accummulation of
infinitely many oscillatory components. Roughly speaking, the construction starts with
a \emph{subsolution}, for which the density of a physical quantity and the corresponding flux satisfy a linear equation, where the strict
constitutive relation is replaced by a kind of convex relaxation. The solution are then obtained by modulating on a given subsolution a family of oscillatory increments that successively ``improve'' the approximate constitutive relation. Being constructed by accumulating oscillations, these weak solutions are likely to violate the admissibility criteria based on energy dissipation, in particular the standard mechanical energy balance. On the other hand, however, the method
produces (non-smooth) initial data, for which the energy inequality and/or similar admissibility criteria do hold. As a matter of fact,
the kinetic energy of these solutions is fully controlled in the process of construction and could be fixed as an arbitrary function of time, see Chiodaroli \cite{Chiod}, Chiodaroli, DeLellis, and Kreml \cite{ChiDelKre}, DeLellis and Sz\' ekelyhidi \cite{DelSze3}.

In 1973, Dafermos \cite{Dafer} proposed an admissibility criterion based on \emph{maximal dissipation} - the physical solutions are those for which the entropy
is produced at the highest possible rate. Several variants of this principle appeared in the literature, see e.g. Krej\v c\' \i \
and Stra\v skraba \cite{KreStr}, and the unique physical solutions were identified in a number of cases, in particular for problems in the simplified 1-D geometry of the physical space, see Dafermos \cite{Dafer}. Although \emph{maximal dissipation} may not be always the desired property for certain models (see Dafermos \cite{Daf55}), it seems relevant in the context of fluid dynamics.

In this paper, we discuss the implications of the principle of \emph{maximal dissipation} in the context of the method of convex integration applied to the
barotropic Euler system. First, extending the results of Chiodaroli \cite{Chiod}, we identify a rich family of initial data for which the method of convex integration yields infinitely many \emph{global-in-time} entropy admissible solutions. Next, we show that neither of these solutions meets the principle of \emph{maximal dissipation}. The last observation can be interpreted as an argument in favour of \emph{maximal dissipation} as a suitable admissibility criterion in mathematical fluid dynamics.

\subsection{Euler system for a barotropic compressible fluid}

Ignoring the thermal effects we introduce the \emph{Euler system} in the form
\bFormula{i1}
\partial_t \vr + \Div (\vr \vu) = 0,
\eF
\bFormula{i2}
\partial_t (\vr \vu) + \Div (\vr \vu \otimes \vu) + \Grad p(\vr) = 0,
\eF
where $\vr = \vr(t,x)$ is the mass density, $\vu = \vu(t,x)$ the velocity, and $p = p(\vr)$ the pressure, describing the time evolution of a compressible
inviscid fluid. To avoid problems connected with the presence of a kinematic boundary, we adopt the standard simplification assuming that the
motion is space-periodic, with the underlying physical domain $\Omega$ identified with the flat torus
\bFormula{i3}
\Omega = \left( [-1,1]|_{\{ -1; 1 \}} \right)^N,\ N=2,3
\eF
The problem (\ref{i1} - \ref{i3}) is supplemented with the initial conditions
\bFormula{i4}
\vr(0, \cdot) = \vr_0, \ \vu(0, \cdot) = \vu_0, \ \vr_0 > 0 \ \mbox{in}\ \Omega.
\eF

\bRemark{ee1}

We have deliberately excluded the case $N=1$, where the method of convex integration fails.
\eR

As is well known, solutions to the system (\ref{i1}), (\ref{i2}) may develop singularities (shocks) in a finite time no matter how smooth and ``small'' the initial data are. Thus if we still believe the system describes the fluid motion also for large times, a kind of
 generalized (weak) solutions must be considered.
We say that $[\vr, \vu]$ is a \emph{weak solution} to the problem (\ref{i1} - \ref{i4}) in $[0,T] \times \Omega$ if:
\bFormula{i5}
\intO{ \Big( \vr(\tau, \cdot) \varphi (\tau, \cdot) - \vr_0 \varphi(0, \cdot) \Big) } =
\int_0^\tau \intO{ \Big( \vr \partial_t \varphi + \vr \vu \cdot \Grad \varphi \Big) } \ \dt
\eF
for any $\tau \in [0,T]$, and any $\varphi \in C^\infty([0,T] \times \Omega)$;
\bFormula{i6}
\intO{ \Big( (\vr \vu)(\tau, \cdot) \cdot \varphi (\tau, \cdot) - \vr_0 \vu_0 \cdot \varphi(0, \cdot) \Big) }
= \int_0^\tau \intO{ \Big( \vr \vu \cdot \partial_t \varphi + \vr \vu \otimes \vu : \Grad \varphi + p(\vr) \Div \varphi \Big) } \ \dt
\eF
for any $\tau \in [0,T]$, and any $\varphi \in C^\infty([0,T] \times \Omega;R^N)$.

The weak solutions considered in this paper are always \emph{bounded} measurable functions, in particular, the quantities
\[
\vr \in C_{\rm weak} ([0,T]; L^1(\Omega)), \ (\vr \vu) \in C_{\rm weak} ([0,T]; L^1(\Omega;R^N))
\]
have well defined instantaneous values $\vr(t, \cdot)$, $(\vr \vu) (t, \cdot)$.

We say that $[\vr, \vu]$ is a weak solution on the time interval $[0, \infty)$ (global-in-time weak solution) if it is a weak solution on any interval $[0,T]$, $T > 0$ finite.

\subsection{Admissible weak solutions}
\label{aw}

The weak solutions of the system (\ref{i1} - \ref{i4}) are, in general, not unique for given initial data. In order to pick up the physically relevant solutions, we recall the mechanical energy equation, multiplying (\ref{i2}) on $\vu$:
\bFormula{i7}
\partial_t \left( \frac{1}{2} \vr |\vu|^2 + P(\vr) \right) + \Div \left[ \left( \frac{1}{2} \vr |\vu|^2 + P(\vr) \right) \vu \right] +
\Div (p(\vr) \vu ) = 0, \ \mbox{with} \ P(\vr) = \vr \int_1^\vr \frac{p(z)}{z^2} \ {\rm d}z.
\eF
We assume that the pressure is an increasing function of the density, specifically,
\bFormula{i7a}
p \in C[0, \infty) \cap C^1(0, \infty), \ p(0) = 0, \ p'(\vr) > 0 \ \mbox{for all}\ \vr > 0.
\eF
Consequently, the potential $P = P(\vr)$ is strictly convex.

Although relation (\ref{i7}) can be deduced from the original system as long as the solution is smooth, it may be violated by some weak solutions of the
same problem. This motivates the introduction of a class of \emph{admissible} weak solutions satisfying the \emph{energy inequality}
\[
\partial_t \left( \frac{1}{2} \vr |\vu|^2 + P(\vr) \right) + \Div \left[ \left( \frac{1}{2} \vr |\vu|^2 + P(\vr) \right) \vu \right] +
\Div (p(\vr) \vu ) \leq 0,
\]
or, more appropriately, its weak form:
\bFormula{i8}
\int_0^T \intO{ \left[ \left( \frac{1}{2} \vr |\vu|^2 + P(\vr) \right) \partial_t \varphi +  \left( \frac{1}{2} \vr |\vu|^2 + P(\vr) \right) \vu
\cdot \Grad \varphi + p(\vr) \vu \cdot \Grad \varphi \right] } \ \dt
\eF
\[
+
\intO{ \left( \frac{1}{2} \vr_0 |\vu_0|^2 + P(\vr_0) \right) \varphi(0,\cdot) } \geq 0
\]
for any $\varphi \in \DC([0,T) \times \Omega)$, $\varphi \geq 0$.

Unlike the state variables $[\vr, \vu]$, the energy
\[
t \mapsto \left( \frac{1}{2} \vr |\vu|^2 + P(\vr) \right)(t, \cdot)
\]
may not be weakly continuous, however, using (\ref{i8}) we can still uniquely identify the limits
\[
E(\tau +), \ E(\tau - ) \in L^\infty(\Omega),
\]
\[
\intO{ E(\tau+) \varphi } = {\rm ess} \lim_{t \to \tau +} \intO{ \left( \frac{1}{2} \vr |\vu|^2 + P(\vr) \right)(t, \cdot) \varphi }
\ \mbox{for}\ \varphi \in C^\infty (\Omega), \ \tau \in [0,T),
\]
and, similarly,
\[
\intO{ E(\tau-) \varphi } = {\rm ess} \lim_{t \to \tau -} \intO{ \left( \frac{1}{2} \vr |\vu|^2 + P(\vr) \right)(t, \cdot) \varphi }
\ \mbox{for}
\ \varphi \in C^\infty (\Omega), \ \tau \in (0,T].
\]
It is easy to check that the mapping $\tau \mapsto E(\tau)$ is a sum of a weakly continuous and a monotone (non-increasing) component, in particular,
\[
E(\tau+) \leq E(\tau-) \ \mbox{for all}\ \tau \in (0,T),
\]
where the identity holds with a possible exception of at most countable set of times. Moreover,
\[
\intO{ E(\tau_2-) } \leq \intO{ E(\tau_1+) } \ \mbox{whenever}\ 0 \leq \tau_1 < \tau_2 \leq T.
\]

\bRemark{aw1}

Since the density $\vr$ and the momentum $(\vr \vu)$ are weakly continuous in $t$, the energy
\[
t \mapsto \frac{1}{2} \vr |\vu|^2 + P(\vr) = \frac{1}{2} \frac{ |\vr \vu|^2 }{\vr} + P(\vr)
\]
is weakly lower semi-continuous, in particular,
\bFormula{i9}
\frac{1}{2} \vr |\vu|^2 + P(\vr) (\tau, \cdot) \leq E(\tau +) \ \mbox{for all}\ \tau \in (0,T),
\eF
and
\[
\frac{1}{2} \vr |\vu|^2 + P(\vr) (\tau, \cdot) = E(\tau-) = E(\tau +) \ \mbox{for a.a.}\ \tau \in (0,T).
\]

Moreover, it follows from the energy inequality (\ref{i8}) that
\[
E(0+) = \frac{1}{2} \vr_0 |\vu_0|^2 + P(\vr_0)
\]

\eR

\subsection{Principle of maximal dissipation}

Adapting slightly the original definition of Dafermos \cite{Dafer},
we say that an admissible weak solution $[\vr, \vu]$ of the Euler system (\ref{i1} - \ref{i4}) satisfies the principle of \emph{maximal dissipation}
if the following holds:

\medskip

{\it Let $\tau \in [0,T)$ and let $[\tilde \vr, \tilde \vu]$ be another weak solution of (\ref{i1} - \ref{i4}), defined in $[0, \tilde T]$,
$\tau < \tilde T \leq T$, such that
\[
\vr = \tilde \vr, \ \vr \vu = \tilde \vr \tilde \vu \ \mbox{in} \ [0,\tau] \times \Omega.
\]

Then there exists a sequence $\{ \tau_n \}_{n=1}^\infty$, $\tau_n > \tau$, $\tau_n \to \tau$ such that
\[
\intO{ \tilde E (\tau_n +) } \geq \intO{ E(\tau_n +) } \ \mbox{for all}\ n = 1,2,\dots,
\]
where $E$, $\tilde E$ is the mechanical energy associated with $u$, $\tilde u$, respectively.

}

\medskip

In other words, the admissible solutions that comply with the principle of \emph{maximal dissipation} loose their mechanical energy at the highest possible rate.

\subsection{Main results}

DeLellis and Sz\' ekelyhidi \cite[Theorem 2]{DelSze3} showed the existence of initial data $[\vr_0, \vu_0]$ for which the problem (\ref{i1}-\ref{i4}) possesses infinitely many global-in-time admissible weak solutions. As a matter of fact, they take $\vr \equiv 1$ and look for solutions of the \emph{incompressible} Euler system with constant
pressure. Chiodaroli \cite{Chiod} used a refined but still similar approach, with the ansatz $\vr = \vr_0$, where $\vr_0$ is a given (smooth) function, and where the solutions have constant (in time) density and the pressure $p(\vr_0)$. In this case, the mechanical energy is changing with time and the energy inequality (\ref{i8}) is valid only on a possibly short interval. As a result, the weak solutions constructed by Chiodaroli \cite{Chiod} are admissible only locally in time.

In this paper we show the existence of global-in-time admissible solutions for any smooth initial distribution of the density $\vr_0$ with sufficiently small
$\Grad \vr_0$, meaning $\vr_0$ exhibits small oscillations around a positive constant state. The precise statement of our result reads:

\bTheorem{i1}

Let the pressure $p = p(\vr)$ satisfy the hypothesis (\ref{i7a}). Let $\vr_0 \in C^1(\Omega)$ be given,
\[
0 < \underline{\vr} \leq \vr_0 (x) \leq \Ov{\vr} \ \mbox{for all}\ x \in \Omega.
\]

Then there exists $\ep > 0$, depending only on the bounds $\underline{\vr}, \Ov{\vr}$, and the structural properties of $p$, and an initial distribution of
the velocity
\[
\vu_0 \in L^\infty (\Omega; R^N)
\]
such that
the problem (\ref{i1} - \ref{i4}) admits infinitely many
admissible weak solutions $[\vr, \vu]$ in $[0, \infty) \times \Omega$ whenever
\[
\sup_{x \in \Omega} |\Grad \vr_0 (x) | < \ep.
\]

\eT

The main idea behind the proof of Theorem \ref{Ti1} is a simple observation, already exploited in \cite{ChiFeiKre}, that the time evolution of the density
$\vr$ is driven by the \emph{acoustic} component of the velocity field, namely
\[
\partial_t \vr = - \Delta \Psi,
\]
where $\Grad \Psi$ is the gradient part in the Helmholtz decomposition of the momentum $\vr \vu$. Thus prescribing {\it a priori} $\vr = \vr(t, \cdot)$, with the associated
acoustic potential $\Psi$, we can construct the desired global-in-time solutions applying the technique of convex integration only to the solenoidal
component of the velocity (momentum).

Next, we turn attention to the way how the solutions are constructed. Roughly speaking, the solenoidal part $\vc{v}$ of the momentum belongs to
the $C_{\rm weak}([0,T], L^2(\Omega;R^N))$-closure of the set $X_{0,e}$ of \emph{subsolutions},
\[
X_{0,e}[0,T] = \left\{ \vc{v} \in C_{\rm weak}([0,T], L^2(\Omega;R^N)) \ \Big| \ \vc{v}(0, \cdot) = \vc{v}_0, \ \vc{v}(T, \cdot) = \vc{v}_T ; \right.
\]
\[
\vc{v} \in C^1((0,T) \times \Omega;R^N), \ \partial_t \vc{v} + \Div \tn{U} = 0 \ \mbox{for a certain}\ \tn{U} \in C^1((0,T) \times \Omega; R^{N \times N}_{\rm sym}),
\]
\[
\left.
\frac{1}{2 \vr} |\vc{v}(t,x) + \Grad \Psi (t,x) |^2 \leq H\Big(\vr, \Psi, \vc{v}(t,x), \tn{U}(t,x) \Big) < e(t,x) \ \mbox{in}\ (0,T) \times \Omega \right\},
\]
where $H = H(\cdot, \vc{v}, \tn{U})$ is a suitable convex function, and $e \in C([0,T] \times \Omega)$ is a given ``kinetic energy'', see Section \ref{p} for details.

The key ingredient in the construction of admissible solutions is the following result:

\bTheorem{i2}
Suppose that $\vc{v} \in X_{0,e}[0,T]$.

Then for any $\tau \in (0,T)$ and any $\ep > 0$ there exists $\Ov{\tau} \in (0,T)$,
\[
|\tau - \Ov{\tau} | < \ep,
\]
and $\vc{w} \in X_{0,e}[\Ov{\tau}, T]$
satisfying $\vc{w}(\Ov{\tau}, \cdot) = \vc{w}_{\Ov{\tau}}$, $\vc{w}(T, \cdot) = \vc{v}_T$,
\bFormula{ccol}
\frac{1}{2} | \vc{w}_{\Ov{\tau}} + \Grad \Psi  (\Ov{\tau}, \cdot) |^2 = e(\Ov{\tau}, \cdot) ,\ \vc{w}(t, \cdot) = \vc{v}(t, \cdot) \ \mbox{in a left neighbourhood of}\ T.
\eF

\eT

As we shall see in Section \ref{a},
the function $\vc{w}$ obtained in Theorem \ref{Ti2} is used as a subsolution in the process of construction of global-in-time admissible solutions to the
problem (\ref{i1} - \ref{i4}). In particular, we immediately deduce the following conclusion (see Section \ref{m} for details and further discussion):

\bCorollary{i1}

The admissible solutions of the problem (\ref{i1} - \ref{i4}) constructed in the proof of Theorem \ref{Ti1} by the method of convex integration
\emph{do not} comply with the principle of maximal dissipation.

\eC

The general outline of the paper is as follows. In Section \ref{p}, we reformulate the problem in the form suitable for a direct application of the method
of convex integration. Section \ref{o} contains a variant of a fundamental lemma on the construction of oscillatory solutions. In Section \ref{a}, we
first prove Theorem \ref{Ti2} and then
complete the proof of Theorem \ref{Ti1}. Finally, Section \ref{m} reveals how the solutions constructed in the
course of the proof of Theorem \ref{Ti1} violate the principle of maximal dissipation. What is more, for each such solution it is possible to
construct another one that dissipates more kinetic energy on any subset of $\Omega$ of positive measure.

\section{Preliminaries}
\label{p}

Similarly to \cite{ChiFeiKre}, it is convenient to express the momentum $\vr \vu$ in terms of its Helmholtz projection $\vc{H}$ and the
acoustic potential $\Psi$, namely
\[
\vr \vc{u} = \vc{v} + \Grad \Psi, \ \vc{v} = \vc{H}[\vr \vu], \ \Div \vc{v} = 0.
\]
Accordingly, we rewrite the system (\ref{i1}), (\ref{i2}) in the form
\bFormula{p1}
\partial_t \vr + \Delta \Psi = 0,
\eF
\bFormula{p2}
\partial_t \vc{v} + \Div \left( \frac{ (\vc{v} + \Grad \Psi) \otimes (\vc{v} + \Grad \Psi) }{\vr} \right) +
\Grad \Big( \partial_t \Psi + p(\vr) - \chi \Big) = 0,
\eF
with a suitable spatially homogeneous function $\chi = \chi(t)$ chosen to satisfy (\ref{o1}) below.

\subsection{Time evolution of the density}
\label{te}

For $\vr_0$ as in Theorem \ref{Ti1}, we take
\bFormula{p3-}
\vr (t,x) =
h(t) \vr_0(x) + (1 - h(t)) \tilde{\vr} ,\ t \in [0,T],
\eF
where
\[
\tilde{\vr} = \frac{1}{|\Omega|}
\intO{ \vr_0 }, \ h \in C^\infty[0,T], \ 0 \leq h \leq 1, \ h(t) = 1 \ \mbox{for}\ t \in  [0, \frac{T}{4}],\
h(t) = 0 \ \mbox{for}\  t \in [\frac{3}{4} T, T].
\]

Accordingly, the potential $\Psi$ is the unique solution of the elliptic problem
\bFormula{p3}
- \Delta \Psi = \partial_t \vr = h'(t) \Big( \vr_0 - \tilde{\vr} \Big), \ \intO{ \Psi } = 0.
\eF
Note that
\[
\Psi (t, \cdot) = 0 \ \mbox{for}\ t \in \left[ 0, \frac{T}{4} \right] \cup \left[ \frac{3}{4}T, T \right].
\]

\subsection{Kinetic energy}

For $\vr$, $\Psi$ satisfying (\ref{p3-}), (\ref{p3}), we introduce the ``kinetic energy''
\bFormula{p4}
e(t,x) = \chi(t) - \frac{N}{2} \partial_t \Psi(t,x) - \frac{N}{2} p(\vr(t,x))
\eF
\[
=  \chi(t) + h''(t) \frac{N}{2}
\Delta^{-1}[ \vr_0 - \Ov{\vr} ]  - \frac{N}{2} p \Big( h(t) \vr_0 (x) + (1 - h(t)) \tilde{\vr} \Big)
,\ e \in C^1([0,T] \times \Omega).
\]

Now, following DeLellis and Sz\' ekelyhidi \cite{DelSze3}, we introduce the necessary tools to apply the machinery of convex integration. Let
$\lambda_{\rm max}[ \tn{A}]$ denote the maximal eigenvalue of a symmetric matrix $\tn{A} \in R^{N \times N}_{\rm sym}$. We define the set
of subsolutions
\[
X_{0,e}[0,T] = \left\{ \vc{v} \in C_{\rm weak}([0,T]; L^2(\Omega;R^N) \ \Big| \ \vc{v} (0, \cdot) = \vc{v}_0, \ \vc{v}(T, \cdot) = \vc{v}_T ,
\ \Div \vc{v} = 0, \right.
\]
\[
\vc{v} \in C^1((0,T) \times \Omega; R^N), \ \partial_t \vc{v} + \Div \tn{U} = 0 \ \mbox{for a certain}\
\tn{U} \in C^1((0,T) \times \Omega; R^{N \times N}_{{\rm sym},0}),
\]
\[
\left. \frac{N}{2} \lambda_{\rm max} \left[ \frac{ (\vc{v} + \Grad \Psi) \otimes (\vc{v} + \Grad \Psi) }{\vr} - \tn{U} \right] < e
\ \mbox{in}\ (0,T) \times \Omega \right\},
\]
where $e$ is given by (\ref{p4}).
Note that (see \cite{DelSze3})
\[
\frac{N}{2} \lambda_{\rm max} \left[ \frac{ (\vc{v} + \Grad \Psi) \otimes (\vc{v} + \Grad \Psi) }{\vr} - \tn{U} \right] \geq \frac{1}{2}
\frac{ |\vc{v} + \Grad \Psi |^2 }{\vr},
\]
where equality holds only if
\[
\tn{U} = \frac{ (\vc{v} + \Grad \Psi) \otimes (\vc{v} + \Grad \Psi) }{\vr} - \frac{1}{N}
\frac{ |\vc{v} + \Grad \Psi |^2 }{\vr} \tn{I}.
\]
The functional $\lambda_{\rm max}$ can be viewed as a norm on the space $R^{N \times N}_{\rm sym,0}$ of symmetric \emph{traceless} tensors.

Finally, we introduce a functional
\[
I [\vc{v}] = \int_0^T \intO{ \left( \frac{1}{2} \frac{ |\vc{v} + \Grad \Psi |^2 }{\vr} - e \right) } \ \dt.
\]
The functional $I$ measures the distance between a subsolution and the topological boundary of the set $X_{0,e}$. As we shall see below, the weak {solutions}
claimed in Theorem \ref{Ti1} correspond to the points of continuity of $I$ on the closure of $X_{0,e}$. The crucial result providing such a conclusion is the oscillatory lemma discussed in the next section.

\section{Oscillatory lemma and convex integration}
\label{o}

We report here a crucial result which can be viewed as a variable-coefficients version of \cite[Proposition 3]{DelSze3}.

\bLemma{o1}

Suppose that $\vc{v} \in X_{0,e}[0,T]$ and that $0 \leq \tau_1 < \tau_2 \leq T$ are given.

Then there exist sequences
\[
\{ \vc{w}_n \}_{n=1}^\infty \subset \DC((\tau_1, \tau_2) \times \Omega; R^N),\
\{ \tn{U}_n \}_{n=1}^\infty \subset \DC((\tau_1, \tau_2) \times \Omega; R^N)
\]
such that the functions $\vc{v} + \vc{w}_n$ belong to $X_{0,e}[0,T]$, with the associated tensor fields $\tn{U} + \tn{U}_n$,
\[
\vc{w}_n \to 0 \ \mbox{in}\ C_{\rm weak}([0,T]; L^2(\Omega;R^N)),
\]
and
\bFormula{P1Da}
\liminf_{n \to \infty} \| \vc{w}_n \|^2_{L^2((0,T) \times \Omega)} \geq \Lambda \int_{\tau_1}^{\tau_2} \intO{ \left(
e - \frac{1}{2} \frac{|\vc{v} + \Grad \Psi|^2}{\vr} \right)^2 } \ \dt, \ \Lambda > 0,
\eF
where the constant $\Lambda$ depends only on $N$ and the norm of the quantities $\vr, \vr^{-1}, \Grad \Psi, e$ in $L^\infty((0,T) \times \Omega)$.

\eL

We refer to \cite[Lemma 3.2]{ChiFeiKre} for a complete proof of Lemma \ref{Lo1}.

With Lemma \ref{Lo1} at hand, we can construct solutions of the problem (\ref{i1} - \ref{i4}) in $[0,T] \times \Omega$ repeating step by step the
arguments of DeLellis and Sh\' ekelyhidi \cite{DelSze3} (see also \cite[Section 3]{ChiFeiKre}).

\medskip

{\bf Step 1.} Taking $\vc{v}_0 = \vc{v}_T = 0$, we fix the function $\chi \in C^1[0,T]$ in (\ref{p4}) in such a way that
\[
\vc{v} = 0 \ \mbox{with}\ \tn{U} = 0 \ \mbox{belong to} \ X_{0,e}[0,T],
\]
specifically, we need
\bFormula{o1}
\chi(t) > \frac{N}{2} \left( \partial_t \Psi (t,x) + p(\vr)(t,x) + \lambda_{\rm max} \left[ \frac{ \Grad \Psi (t,x) \otimes \Grad \Psi (t,x) }{\vr (t,x) } \right] \right)
\ \mbox{for all}\ (t,x) \in [0,T] \times \Omega.
\eF
In particular, the set of subsolutions $X_{0,e}[0,T]$ is non-empty.

\medskip

{\bf Step 2.} Using the oscillatory lemma (Lemma \ref{Lo1}) we show that solutions of (\ref{p1}), (\ref{p2}) can be identified with the points
of continuity of the functional $I$ on the closure of $X_{0,e}[0,T]$ in $C_{\rm weak}([0,T]; L^2(\Omega;R^N))$. In particular, these solutions
satisfy
\bFormula{o2}
\vc{v}(0, \cdot) = \vc{v}(T, \cdot) = 0, \ \vr = \vr_0 \ \mbox{in} \ [0, \frac{T}{4}],\ \vr = \tilde{\vr} \ \mbox{in} \ [ \frac{3}{4}T,T],
\eF
\bFormula{o3}
\frac{1}{2} \frac{ |\vc{v} + \Grad \Psi |^2 }{\vr} = e =  \chi(t) - \frac{N}{2} \partial_t \Psi(t,x) - \frac{N}{2} p(\vr(t,x))
\ \mbox{a.a. in}\ (0,T) \times \Omega,
\eF
\[
\mbox{where}\ \Psi = 0 \ \mbox{in}\ [0, \frac{T}{4}] \cup [\frac{3T}{4},T].
\]

\medskip

{\bf Step 3.} Exactly as in
\cite[Section 3]{ChiFeiKre} we conclude that there are infinitely many solutions to the problem (\ref{p1}), (\ref{p2}) - the points of continuity of $I$ -  satisfying (\ref{o2}),
(\ref{o3}) provided $\chi$ is chosen as in (\ref{o1}).

\bRemark{oo1}

The solutions constructed in the above are obviously not admissible, as they start from the initial state $\vc{v}_0 = 0$ and
\emph{produce} mechanical energy. They can be viewed as analogues of the solutions to the incompressible Euler system obtained by
Scheffer \cite{Scheff} and Shnirelman \cite{Shn}.

\eR

In order to obtain \emph{admissible weak solutions}, we have to construct subsolutions, with the associated mechanical (kinetic) energy left-continuous at the point $t = 0$. This will be done in the next section.

\section{Infinitely many admissible solutions}
\label{a}

Our goal in this section is to complete the proof of Theorem \ref{Ti1}. To this end, we first identify a class of suitable initial data
for which we can construct subsolutions with the associated mechanical energy continuous at $t = 0$. This is basically the claim of Theorem
\ref{Ti2} so we we start by its proof, which can be seen as an extension of \cite[Proposition 5]{DelSze3} to the case of non-constant energy and variable
coefficients. The proof is quite illustrative since it shows how the machinery of convex integration works, accumulating small oscillatory
increments to a given subsolution.

\subsection{Proof of Theorem \ref{Ti2}}

For a fixed energy function $e \in C^1([0,T] \times \Omega)$, we identify a suitable class of initial velocities $\vc{v}_0$, for which
the kinetic energy is right-continuous at $t = 0$ (cf. (\ref{ccol}). In accordance with the hypotheses of Theorem \ref{Ti2}, we suppose that the set of subsolutions $X_{0,e}[0,T]$ contains at least
one element $\vc{v}$.

\subsubsection{Oscillatory sequence}

In order to construct a suitable subsolution, we apply successively the oscillatory lemma (Lemma \ref{Lo1}). More specifically, the
function $\vc{w}$, the existence of which
is claimed in Theorem \ref{Ti2}, is obtained as a limit of a sequence $\{ \vc{w}_k \}_{k=1}^\infty \subset X_{0,e}[0,T]$,
\[
\vc{w}_k \to \vc{w} \ \mbox{in}\ C_{\rm weak} ([0,T]; L^2(\Omega;R^N)),
\]
where we take
\[
\vc{w}_0 = \vc{v} \equiv \vc{v}_0 ,
\]
and fix
\[
\tau_0 = \tau,\ \ep_0 = \ep.
\]
The functions $\vc{w}_k$ are defined recursively by the following procedure:

\begin{enumerate}

\item the increment $\vc{w}_k - \vc{w}_{k-1}$ is compactly supported in a small neighborhood  of the point $\tau_{k-1}$, specifically,

\bFormula{AA3-}
\vc{w}_{k} \in X_{0,e}[0,T],\ {\rm supp}[\vc{w}_k - \vc{w}_{k-1}] \subset (\tau_{k - 1} - \ep_k, \tau_{k - 1} + \ep_k),\
 \mbox{where}\ 0< \ep_k < \frac{\ep_{k-1}}{2};
\eF

\item $\vc{w}_{k} - \vc{w}_{k-1}$ is small in the topology of the space $C_{\rm weak}([0,T]; L^2(\Omega;R^N))$, we take

\bFormula{AA2-}
d(\vc{w}_k, \vc{w}_{k-1}) < \frac{1}{2^k}, \ \sup_{t \in (0,T)} \left| \intO{ \frac{1}{\vr} (\vc{w}_k - \vc{w}_{k-1}) \cdot \vc{w}_{m} } \right| < \frac{1}{2^k}
\ \mbox{for all}\ m=0,\dots, k-1,
\eF
for $k=1, \dots$, where
the symbol $d$ denotes the metric induced by the weak topology on bounded sets of $L^2(\Omega;R^N)$;

\item
the sequence $\{ \vc{w}_k \}_{k=0}^\infty$ is oscillatory around a point $\tau_k$,
\[
\tau_{k} \in (\tau_{k-1} - \ep_k,
\tau_{k-1} + \ep_k),
\]
meaning
\bFormula{AA2}
\intO{ \frac{1}{2} \frac{|\vc{w}_{k} + \Grad \Psi|^2}{\vr} (\tau_k, \cdot) }  \geq
\intO{ \frac{1}{2} \frac{|\vc{w}_{k-1} + \Grad \Psi|^2}{\vr} (t, \cdot) }  + \frac{\lambda}{\ep^2_k} \alpha_k^2
\eF
\[
\geq \intO{ \frac{1}{2} \frac{|\vc{w}_{k-1} + \Grad \Psi |^2}{\vr} (\tau_{k-1}, \cdot) }  + \frac{\lambda}{2 \ep^2_k} \alpha_k^2 \ \mbox{for all}\ t \in (\tau_{k-1} - \ep_k, \tau_{k - 1} + \ep_k),
\]
where
\[
\alpha_k = \int_{\tau_{k-1} - \ep_k}^{\tau_{k-1} + \ep_k} \intO{ \left( e - \frac{1}{2} \frac{ |\vc{w}_{k-1} + \Grad \Psi|^2 } {\vr} \right) } \ \dt > 0,
\]
and $\lambda > 0$ is constant independent of $k$.

\end{enumerate}

With $\vc{w}_0, \dots, \vc{w}_{k - 1}$ already constructed, our goal is to find $\vc{w}_{k}$ enjoying the properties (\ref{AA3-} - \ref{AA2}).
To this end, we compute
\[
\alpha_k = \int_{\tau_{k-1} - \ep_k}^{\tau_{k-1} + \ep_k} \intO{ \left( e - \frac{1}{2} \frac{ |\vc{w}_{k-1} + \Grad \Psi|^2 } {\vr} \right) } \ \dt
\ \mbox{for} \ 0 < \ep_k < \frac{\ep_{k-1}}{2}
\]
and observe that
\[
\frac{\alpha_{k}}{2 \ep_{k}} = \frac{1}{2 \ep_{k}} \int_{\tau_{k-1} - \ep_{k}}^{\tau_{k-1} + \ep_{k}} \intO{ \left( e - \frac{1}{2} \frac{ |\vc{w}_{k - 1} +
\Grad \Psi |^2 } {\vr} \right) } \ \dt
\]
\[
\to \intO{ \left( e - \frac{1}{2} \frac{ |\vc{w}_{k-1} + \Grad \Psi|^2 } {\vr} \right)(\tau_{k-1}) } > 0
\ \mbox{for}\ \ep_{k} \to 0
\]
as the function $\vc{w}_{k-1}$ is smooth in $(0,T)$.

Next, we take $\ep_k > 0$ small enough so that
\bFormula{odhad}
\frac{1}{2 \ep_k} \int_{\tau_{k-1} - \ep_{k} }^{ \tau_{k-1} + \ep_{k}}
\intO{ \frac{1}{2} \frac{|\vc{w}_{k-1} + \Grad \Psi |^2}{\vr}  } \ \dt + \frac{\Lambda(\| \vr, \vr^{-1}, \Grad \Psi, e \|_{L^\infty(0,T) \times \Omega)}) }{4 \ep_{k}^2} \alpha_{k}^2
\eF
\[
\geq
\intO{ \frac{1}{2} \frac{|\vc{w}_{k-1} + \Grad \Psi |^2}{\vr}(t,\cdot)  }  + \frac{\Lambda(\| \vr, \vr^{-1}, \Grad \Psi, e \|_{L^\infty(0,T) \times \Omega)})}{8 \ep_{k}^2} \alpha_{k}^2
\]
\[
\geq \intO{ \frac{1}{2} \frac{|\vc{w}_{k-1} + \Grad \Psi|^2}{\vr}(\tau_{k-1},\cdot)  }  + \frac{\Lambda(\| \vr, \vr^{-1}, \Grad \Psi, e \|_{L^\infty(0,T) \times \Omega)})}{16 \ep_{k}^2} \alpha_{k}^2
\]
\[
\mbox{for all}\ t \in (\tau_{k-1} - \ep_k, \tau_{k-1} + \ep_k),
\]
where $\Lambda(\| \vr, \vr^{-1}, \Grad \Psi, e \|_{L^\infty(0,T) \times \Omega)}) > 0$ is the universal constant from Lemma \ref{Lo1}.

Applying Lemma \ref{Lo1} we can construct a function $\vc{w}_k \in X_{0,e}$
such that
\[
{\rm supp}[ \vc{w}_{k} - \vc{w}_{k-1} ] \subset (\tau_{k-1} - \ep_{k}, \tau_{k-1} + \ep_{k} ),
\]
\bFormula{aA1}
d(\vc{w}_{k}, \vc{w}_{k-1} ) < \frac{1}{2^{k}}, \ \sup_{t \in (0,T)} \left| \intO{ \frac{1}{\vr} (\vc{w}_{k} - \vc{w}_{k-1}) \cdot \vc{w}_{m} } \right| < \frac{1}{2^k},\ m=0,\dots, k-1
\eF
and
\bFormula{odhad2}
\int_{\tau_{k-1} - \ep_{k} }^{ \tau_{k-1} + \ep_{k}} \intO{ \frac{1}{2} \frac{|\vc{w}_{k} + \Grad \Psi|^2}{\vr}  } \ \dt
\geq
\int_{\tau_{k-1} - \ep_{k} }^{ \tau_{k-1} + \ep_{k}}
\intO{ \frac{1}{2} \frac{|\vc{w}_{k-1} + \Grad \Psi|^2}{\vr}  } \ \dt + \frac{\Lambda}{2 \ep_{k}} \alpha_{k}^2,
\eF
where the last integral in (\ref{P1Da}) has been estimated from below by  Jensen's inequality.
The relations (\ref{odhad}), (\ref{odhad2}) yield (\ref{AA2}) for some $\tau_k \in (\tau_{k-1} - \ep_k, \tau_{k-1} + \ep_k)$, and with
$\lambda = \Lambda / 16$.

\subsubsection{A suitable subsolution}

It follows from (\ref{AA2-}) that there is $\vc{w}$ such that
\bFormula{conv}
\vc{w}_k \to \vc{w} \ \mbox{in} \ C_{\rm weak}([0,T]; L^2(\Omega; R^N)).
\eF
Moreover, by virtue of (\ref{AA3-}),
\[
\tau_k \to \Ov{\tau} \in (0,T) ,\ | \Ov{\tau} - \tau | < \ep.
\]
Finally, for any $\delta > 0$, there is $k = k_0(\delta)$ such that
\bFormula{AA5}
\vc{w}(t, \cdot) = \vc{w}_{k} (t, \cdot) = \vc{w}_{k_0}(t, \cdot) \ \mbox{for all} \ t \in (0, \Ov{\tau} - \delta) \cup (\Ov{\tau} + \delta, T),
k \geq k_0.
\eF
In particular, we deduce from (\ref{AA5}) that
\[
\vc{w} \in X_{0,e}[\Ov{\tau},T] \ \mbox{with}\
\vc{w}(\Ov{\tau}, \cdot) \in L^\infty(\Omega;R^N), \ \mbox{and}\ \vc{w} \equiv \vc{v}, \ \tn{U}_{\vc{w}} \equiv 0 \ \mbox{in a (left) neighborhood of}\ T.
\]

\subsubsection{Continuity of the initial energy}

Our ultimate goal is to show that the subsolution $\vc{w}$ is continuous from the left at $\Ov{\tau}$ in the strong topology of $L^2(\Omega;R^N)$,
which is equivalent to (\ref{ccol}).
To this end, we first observe that
(\ref{AA2}) implies
\bFormula{AA6-}
\intO{ \frac{1}{2} \frac{ |\vc{w}_{k-1} + \Grad \Psi |^2 }{\vr} (t, \cdot) } \nearrow Z \ \mbox{uniformly for}\ t \in (\tau_{k-1} - \ep_k, \tau_{k-1} + \ep_k),
\eF
therefore
\bFormula{AA6}
\frac{\alpha_k}{\ep_k} = \frac{1}{\ep_k} \int_{\tau_{k-1} - \ep_k}^{\tau_{k-1} + \ep_k} \intO{ \left( e - \frac{1}{2} \frac{ |\vc{w}_{k-1} + \Grad \Psi|^2 } {\vr} \right) } \ \dt \to 0;
\eF
whence, finally,
\bFormula{AA7}
\intO{ \frac{1}{2} \frac{ |\vc{w}_k + \Grad \Psi|^2 }{\vr} (\Ov{\tau} , \cdot) } \nearrow \intO{ e(\Ov{\tau}, \cdot) }.
\eF

Combining the relations (\ref{AA7}) with (\ref{AA2-}), (\ref{conv}) we may infer that
\[
\vc{w}_k (\Ov{\tau}, \cdot) \to \vc{w}(\Ov{\tau}, \cdot) \ \mbox{in}\ L^2(\Omega; R^N)
\]
which yields (\ref{ccol}). Indeed we have
\[
\intO{ \frac{1}{\vr} | \vc{w}_n - \vc{w}_m |^2 (\Ov{\tau}, \cdot) }
\]
\[
=
\intO{ \frac{1}{\vr} |\vc{w}_n |^2 (\Ov{\tau}, \cdot) } - \intO{ \frac{1}{\vr} |\vc{w}_m |^2 (\Ov{\tau}, \cdot) }
- 2 \intO{ \frac{1}{\vr} \left( \vc{w}_n - \vc{w}_m \right) \cdot \vc{w}_m (\Ov{\tau}, \cdot) }
\ \mbox{for all} \ n > m,
\]
where, by virtue of (\ref{AA2-}),
\[
\intO{ \frac{1}{\vr} \left( \vc{w}_n - \vc{w}_m \right) \cdot \vc{w}_m (\Ov{\tau}, \cdot) }
= \sum_{ k = 0}^{n-m-1} \intO{ \frac{1}{\vr} \left( \vc{w}_{k+1} - \vc{w}_k \right) \cdot \vc{w}_m (\Ov{\tau}, \cdot) }
\to 0 \ \mbox{for}\ m \to \infty.
\]

We have proved Theorem \ref{Ti2}.

\subsection{Proof of Theorem \ref{Ti1}}

With Theorem \ref{Ti2} at hand, we may use the machinery of convex integration to produce admissible weak solutions to the
compressible Euler system, meaning the weak solution that dissipate mechanical energy. To this end, we first refine our
requirements concerning the function $e$ appearing in the definition of $X_{0,e}$.

\subsubsection{Mechanical energy revisited}

The weak solutions constructed in Section \ref{o} as limits of subsolutions in $X_{0,e}[0,T]$ have the energy defined
(\ref{o3}), with $\chi$ satisfying (\ref{o1}).

Accordingly, the
energy inequality (\ref{i8}) takes the form
\bFormula{o4}
\partial_t \chi (t) - \frac{N}{2} \partial^2_{t,t} \Psi (t,x) - \frac{N}{2} \partial_t p(\vr)(t,x) + \partial_t P(\vr)(t,x) =
\partial_t \left[ e + P(\vr) \right]
\eF
\[
\leq - \Div \left[ \left( e + P(\vr) + p(\vr) \right) \frac{\vc{v} + \Grad \Psi}{\vr} \right]
\]
\[
=
- \chi(t) \Div \left[ \frac{\vc{v} + \Grad \Psi}{\vr} \right] + \Div \left[ \left( \frac{N}{2} \partial_t \Psi + \left( \frac{N}{2} - 1\right) p(\vr) + P(\vr) \right)
\frac{\vc{v} + \Grad \Psi}{\vr} \right].
\]
We note that, in view of (\ref{o2}), (\ref{o3}), the relation (\ref{o4}) reduces to
\bFormula{o5}
\partial_t \chi (t) \leq 0 \ \mbox{as soon as}\ \frac{3}{4} T \leq t \leq T.
\eF

On the other hand, since $\vr$, $\Psi$ were fixed by (\ref{p3-}), (\ref{p3}), relation (\ref{o4}) follows as soon as $\chi$ satisfies
\bFormula{oo5}
\partial_t \chi \leq - \chi(t) \left[ C_1(\vr, \Psi) + \Div \left( \frac{\vc{v}}{\vr} \right) \right] + \vc{v} \cdot \vc{C}_2(\vr, \Psi) - C_3(\vr, \Psi) ,\ t \in (0,T).
\eF
Moreover, the relation (\ref{o3}) yields
\[
| \vc{v}(t,x) | \leq C_5 (\vr, \Psi ) \left( 1 + \sqrt{ \chi (t) } \right).
\]

Thus, seeing that $\Div \vc{v} = 0$, we can find a smooth function $\chi$ satisfying both (\ref{o1}) and (\ref{o4}) on the whole time
interval $[0,T]$ as soon as
\[
\sup_{x \in \Omega} |\Grad \vr_0 (x)| < \ep ,\ \mbox{with}\ \ep > 0 \ \mbox{small enough,}
\]
in accordance with the hypotheses of Theorem \ref{Ti1}.

\subsubsection{Adjusting the initial data}

At this stage, we are ready to complete the proof of Theorem \ref{Ti1}. To this end, we fix the functions $\vr$, $\Psi$ as in (\ref{p3-}), (\ref{p3}), and  choose $\chi$ and $e$ satisfying (\ref{o1}), (\ref{o4}). In view of (\ref{o5}), the functions $\chi$ and $e$ can be extended by constant
beyond the time $T$ so that (\ref{o1}), (\ref{o4}) remain valid on the unbounded interval $[0, \infty)$.

At this stage, we take $\tau \in [0, \frac{T}{8})$ and set
$\vc{v}_0 = \vc{w} (\Ov{\tau})$, $\vc{v}_T = 0$, where $\vc{w}$ is the subsolution constructed in Theorem \ref{Ti2} belonging to the set
$X_{0,\ep}[\Ov{\tau}, T]$, $\Ov{\tau} \in [0, \frac{T}{8}]$.

Finally, performing a simple time shift $\Ov{\tau} \to 0$
we define a space of subsolutions $X_{0,e}[0,T]$, where $e$ satisfies  (\ref{o4}), and
\bFormula{ad1}
\vr(0, \cdot) = \vr_0,\ \vr(T, \cdot) = \tilde{\vr},\ \vc{v}_0 = \vc{w} (\Ov{\tau}),\
\frac{1}{2} \frac{|\vc{v}_0 + \Grad \Psi(0, \cdot)|^2}{\vr_0}  = \frac{1}{2} \frac{|\vc{v}_0 |^2}{\vr_0} = e(0, \cdot),\
\vc{v}_T = 0.
\eF

By virtue of Theorem \ref{Ti2}, the space $X_{0,e}[0,T]$ is non-empty; whence we may use
the arguments of convex integration
discussed in Section \ref{o} to construct \emph{infinitely
many} admissible weak solutions to the problem (\ref{p1}), (\ref{p2}), or, equivalently, to the problem (\ref{i1} - \ref{i4})
on the time interval $[0, \Ov{T}]$ satisfying
\[
\vr(0, \cdot) = \vr_0, \vr \vu (0, \cdot) = \vr_0 \vc{v}_0, \ \vr(\Ov{T}, \cdot) = \tilde{\vr}, \ \vu(T, \cdot) = 0
\]
as long as $|\Grad \vr_0|$ is small enough. We recall that these solutions are indeed admissible, as they satisfy the energy inequality
in $(0,T)$ thanks to (\ref{o4}), while the initial energy is correctly adjusted by (\ref{ad1}).

Finally, it is easy to check that each of these solutions, extended to be
\[
\vr (t,x) = \tilde{\vr} , \ \vu(t,x) = 0 \ \mbox{for}\ t > \Ov{T} ,\ x \in \Omega,
\]
is a global-in-time admissible weak solution of the same problem. Thus we have completed the proof of Theorem \ref{Ti1}.

\bRemark{o1}

As a matter of fact, the construction of the density in the proof of Theorem \ref{Ti1} could be modified in order to reach an arbitrary density distribution
$\vr(T)$
at the time $T$ satisfying only natural the natural mass constraint
\[
\intO{ \vr_0} = \intO{ \vr(T, \cdot) }.
\]
In particular, we have solved a kind of \emph{control problem}:

\medskip

For a given initial distribution of the density $\vr_0$, find a velocity field $\vu_0$ such that the corresponding solution of the compressible
Euler system reaches the equilibrium state $[\tilde{\vr}, 0]$ at the time $T$, see  Nersisyan \cite{Ner}, and Ervedoza et al. \cite{EGGP} for related results in the viscous case.

\eR

\bRemark{o2}

We have constructed the global-in-time solutions by extending them to coincide with the equilibrium $[\Ov{\vr}, 0]$
beyond the time $T$. Using the technique od DeLellis and Sz\' ekelyhidi \cite{DelSze3} for the incompressible Euler system, we may construct a different
extension (as a matter of fact inifinitely many), namely
\[
\vr = \tilde {\vr}, \ \vu(T, \cdot) = 0, \ |\vu (t, \cdot) | = \alpha > 0 \ \mbox{for all}\ t > T,
\]
where $\alpha$ is chosen small enough so that
\[
E(T-) \geq E(T+) = \intO{ \frac{1}{2} \Ov{\vr} \alpha^2 + P(\tilde{\vr})}.
\]
Obviously, the kinetic (mechanical) energy is not continuous at the point $T$, we have
\[
\intO{ \Ov{\vr} | u(T,\cdot) |^2 + P(\tilde{\vr})  } < \min\{ E(T-), E(T+) \}.
\]

\eR

\section{Maximal dissipation}
\label{m}

The global-in-time solutions with \emph{positive} kinetic energy discussed in Remark \ref{Ro2} obviously \emph{do not} satisfy the principle of maximal dissipation. Clearly, the admissible solutions that coincide with the equilibrium state $[\tilde{\vr}, 0]$ dissipate more energy beyond the time
$T$. As a matter of fact, we can show the following result:

\bProposition{m1}

Let $[\vr, \vu]$ be an admissible weak solution of the problem (\ref{i1} - \ref{i4}) in $[0,T] \times \Omega$ satisfying the principle of maximal dissipation. Suppose that the system (\ref{i1} - \ref{i2}) admits a smooth solution $[\tilde \vr, \tilde \vu]$ in $[\tau_1, \tau_2]$,
$0 \leq \tau_1 < \tau_2 \leq T$ such that
\[
\vr(\tau_1, \cdot) = \tilde \vr(\tau_1, \cdot), \ \vu(\tau_1, \cdot) = \tilde \vu(\tau_1, \cdot).
\]

Then
\[
\vr = \tilde \vr , \ \vu = \tilde \vu \ \mbox{in} \ [\tau_1, \tau_2].
\]

\eP

\bProof

If $\tau_1 = 0$,
the proof reduces to the standard \emph{weak-strong} uniqueness principle (see Dafermos \cite{Daf4}).

Suppose that $\tau_1 > 0$. In order to apply the weak-strong uniqueness result, we have only to show that
\[
\intO{ \left[ \frac{1}{2} \frac{ | (\vr \vu) |^2 }{\vr} + P(\vr) \right] (\tau_1, \cdot) } \geq \intO{ E(\tau_1 +) } =
{\rm ess} \lim_{t \to \tau_1+} \intO{ \left[ \frac{1}{2} \frac{ | (\vr \vu) |^2 }{\vr} + P(\vr) \right] (t, \cdot) }.
\]

Assuming the contrary we would have
\[
\lim_{t \to \tau_1+} \intO{ E(t +) } > \intO{ \left[ \frac{1}{2} \frac{ | (\vr \vu) |^2 }{\vr} + P(\vr) \right] (\tau_1, \cdot) } =
\intO{ \left[ \frac{1}{2} \frac{ | (\tilde \vr \vu) |^2 }{\tilde \vr} + P(\tilde \vr) \right] (t, \cdot) } \ \mbox{for all}\ t \in [\tau_1, \tau_2]
\]
as the smooth solution conserves the energy. This contradicts the principle of maximal dissipation for $[\vr, \vu]$.

\qed

In the light of Proposition \ref{Pm1}, it is obvious that only the solutions obtained in the proof of Theorem \ref{Ti1} that coincide with the
equilibrium after certain time might comply with the principle of maximal dissipation. However, even this is false as we can easily deduce from the process of their construction. Inspecting step by step the proof of Theorem \ref{Ti1} we recall the procedure of convex integration:

\begin{enumerate}
\item fix $\vr$ and the acoustic potential $\Psi$;
\item find suitable $\chi = \chi(t) \in C^1[0,T]$, with the associated kinetic energy $e$,
\[
e =  e(t,x) = \chi(t) - \frac{N}{2} \partial_t \Psi(t,x) - \frac{N}{2} p(\vr(t,x))
\]
so that the weak solutions constructed by the method of convex integration may be admissible;
\item find a suitable initial solenoidal part of the velocity $\vc{v}_0$ satisfying
\[
\frac{1}{2} \frac{|\vc{v}_0 + \Grad \Psi(0, \cdot) |^2}{\vr_0}  = e(0, \cdot) \ \mbox{in}\ \Omega
\]
and such that the set of subsolutions
\[
X_{0,e}[0,T] = \left\{ \vc{v} \in C_{\rm weak}([0,T]; L^2(\Omega;R^N) \ \Big| \ \vc{v} (0, \cdot) = \vc{v}_0, \ \vc{v}(T, \cdot) = 0 ,
\ \Div \vc{v} = 0, \right.
\]
\[
\vc{v} \in C^1((0,T) \times \Omega; R^N), \ \partial_t \vc{v} + \Div \tn{U} = 0 \ \mbox{for a certain}\
\tn{U} \in C^1((0,T) \times \Omega; R^{N \times N}_{{\rm sym},0}),
\]
\[
\left. \frac{N}{2} \lambda_{\rm max} \left[ \frac{ (\vc{v} + \Grad \Psi) \otimes (\vc{v} + \Grad \Psi) }{\vr} - \tn{U} \right] < e
\ \mbox{in}\ (0,T) \times \Omega \right\}
\]
be non-empty;

\item apply the argument of convex integration to deduce the existence of infinitely many admissible weak solutions $[\vr, \vu]$ satisfying
\bFormula{m1}
\frac{1}{2} \vr |\vu|^2 = e \ \mbox{a.a. in}\ (0,T) \times \Omega, \ \vr(0, \cdot) = \vr_0, \vu(0, \cdot) = \vc{v}_0 + \Grad \Psi(0, \cdot).
\eF

\end{enumerate}

As the function $\chi$ (and $e$) is continuous in $t$ (even continuously differentiable), we can choose $\tilde \chi$ with the associated
$\tilde e$ in such a way that
\[
\chi(0) = \tilde \chi(0), \ \tilde \chi(t) < \chi(t) \ \mbox{in} \ (0,T),
\]
with a \emph{non-empty} set of subsolutions $X_{0,\tilde e}[0,T]$. Clearly, the weak solutions obtained from
the space $X_{0, \tilde e}[0,T]$ by the process of convex integration satisfy the same system of equations with the same initial data, however, in accordance with (\ref{m1}), with a strictly larger
rate of dissipation than those obtained from $X_{0, e}[0,T]$. We may therefore infer that solutions constructed in such a way \emph{cannot} comply with the principle of maximal dissipation.

\subsection{Maximal dissipation revisited}

We finish our discussion by introducing a weaker notion of maximal dissipation. Let $[\vr, \vu]$, $[\tilde \vr, \tilde \vu]$ be two admissible weak solutions
of the problem (\ref{i1} - \ref{i4}) defined on $[0,T]$, $[0, \tilde T]$, respectively. We say that
\[
[\vr, \vu]   \succ  [\tilde \vr, \tilde \vu]
\]
if there exists $\tau \in [0, \min\{ T, \tilde T \} )$ such that
\[
\vr(t, \cdot) = \tilde \vr (t, \cdot),\ \vu(t, \cdot) = \tilde \vu (t, \cdot) \ \mbox{for}\ t \in [0, \tau],
\]
and
\[
\tilde E(t+)  > E(t+)
\ \mbox{for any} \ t \ \mbox{in a right neighborhoood of}\ \tau.
\]
Here $h > g$ means $h > g$ a.a. in $\Omega$. If $[\vr, \vu]   \succ  [\tilde \vr, \tilde \vu]$, the
solution $[\vr, \vu]$ dissipates more energy than $[\tilde \vr, \tilde \vu]$
for $t > \tau$, and this inequality is strict on \emph{any} subset of $\Omega$ of positive measure.

It is easy to check that any admissible weak solution $[\vr, \vu]$ satisfying the principle of maximal dissipation is \emph{maximal} with respect to the
relation $\succ$; there is no solution ``greater'' than $[\vr, \vu]$, meaning there is no solution that dissipates more energy. On the other hand,
a direct inspection of the arguments presented at the beginning of this section reveals that the solutions constructed by the method of convex integration \emph{are not} maximal with respect to the relation $\succ$.

\def\cprime{$'$} \def\ocirc#1{\ifmmode\setbox0=\hbox{$#1$}\dimen0=\ht0
  \advance\dimen0 by1pt\rlap{\hbox to\wd0{\hss\raise\dimen0
  \hbox{\hskip.2em$\scriptscriptstyle\circ$}\hss}}#1\else {\accent"17 #1}\fi}


\end{document}